\documentclass[10pt]{article}
\usepackage{mathtools,amssymb,amsmath,verbatim}

\begin{document}

\begin{center}
\textbf{\Large Errata to:  A New Graph over\\ Semi-Direct Products of Groups}\\[8pt]
John K.~McVey\\ Kent State University \footnote{Department of Mathematical Sciences, Kent State University, 233 MSB, 1300 Lefton Esplanade, Kent, Ohio 44242}\\[2pt]\texttt{jmcvey@math.kent.edu}
\end{center}

\smallskip

\begin{abstract}
The goal of the paper ``A new graph over semi-direct products of groups" is to define a graph $\Gamma(G)$ on
a group $G$ when $G$ splits over a normal subgroup.  We demonstrate herein that the graph is ill-defined.  We also
attempt to ascertain causes for the discrepancies.\\
\textbf{MSC: }Primary 20A99,\;05C12
\end{abstract}

\section{Background}
In group theory, a familiar theme is the application of group invariants that can 
be used to demonstrate two groups are not isomorphic.  
In this vein, placing a graph structure on groups and showing corresponding graph 
invariants differ shows the graphs are distinct, hence the underlying groups are 
nonisomorphic. The goal of the paper~\cite{wrong} is such a graph.  There however
is a huge problem with the paper; the graph it presents is ill-defined.  In fact, 
the paper explicitly presents two isomorphic groups then computes
distinct graph invariants thereof.  Hence, the purpose of this paper is to 
demonstrate the graph is in fact ill-defined.  

The graph in question is
described as being on groups $G$ which are split extensions $G=K\rtimes
A$.  Letting $[X;\mathbf r]$ and $[Y;\mathbf s]$ be respectively
presentations for $K$ and $A$ (with $X$ the generators and \textbf r the
relations), the graph $\Gamma(G)$ has vertex set $G$ and (quoting the
paper directly) ``\dots the \emph{edge set} $E$ is obtained by the
following steps:

\medskip

\quad\textbf {(I)} Each of the vertices in this graph must be adjoined to
the vertex~$1_G$ (except~$1_G$ itself since the graph is assumed to be
simple).

\medskip

\quad\textbf{(II)(\textit i)} For any two vertices $w_1 =
x_1^{\varepsilon_1}x_2^{\varepsilon_2}\dots x_m^{\varepsilon_m}$ and $w_2
= y_1^{\delta_1}y_2^{\delta_2}\dots y_n^{\delta_n}$ (where $n\geq 2$,
$\varepsilon_i$ and $\delta_i$ are integers) and for all $x_i,y_j\in X\cup
Y~(1\leq i\leq m, 1\leq j\leq n)$, if $x_i\not=y_j$, then $w_1$
is adjoined to the $w_2$ (shortly, $w_1\sim w_2$).

\medskip

\quad\textbf{~~\;(\textit{ii})} As a consequence of \textbf{(\textit i)},
for any two vertices $w_1 = x_i{}^k$ and $w_2 = x_j{}^t~(1 \leq i, j\leq n,~i\not=j$, 
and $k, t$ are integers), we can directly take $w_1 \sim w_2$. However, to adjoin 
$w_1$ and $w_2$ while $i = j$, it must be $k\not=t$."

\medskip

Clearly, whether or not an edge connects two vertices is strongly tied to
the presentation.  To mitigate this dependence, the article states ``\dots all elements
$z_i~(i=1,2,\dots,k)$ in the generating set $X\cup Y$ of~$G$ will be
formed as $z_i\not=z_1^{\varepsilon_1}z_2^{\varepsilon_2}\dots
z_k^{\varepsilon_k}$, where $k\geq 2$ according to the \emph{Normal Form
Theorem} (\textsl{NFT}) (see [11])." (The article's reference~[11] is our
reference \cite{cohen}.)

\section{Counterexample}

In trying to digest the above definition (wherein $0$ exponents may or may
not be allowed, wherein order of multiplicands is not specified, wherein 
\textbf{(II)(\textit{ii})} claims the adjacency of $x_i{}^k$ and $x_j{}^t$ when 
$i\not=j$ follows from \textbf{(II)(\textit i)} while \textbf{(II)(\textit i)} has a 
hypothesis that there be at least $n\geq 2$ generators in at least one of the vertices, etc.),
this article's author turned to the examples presented in~\cite{wrong} to
see what interpretations were used there; specifically, Examples~2.9
and~2.10.  Example~2.9 is given as the dihedral group $D_8$ viewed as the
split extension of the cyclic $C_4$ being acted upon by~$C_2$ (the action,
of course, being inversion).  Meanwhile, Example~2.10 is the group $G$
defined as the split extension of the Klein~$4$-group $\mathcal V_4$ being
acted upon by $C_2$.  Since ``\dots the homomorphism $\varphi$ will always
be not identity $id_G$ unless stated otherwise," we conclude this is a
nontrivial action, in which case $C_2$ exchanges two of the nontrivial 
members, while fixing the third.  The article~\cite{wrong} calculates the 
respective degree sequences as
$$\operatorname{DS}(\Gamma(D_8)) = \{ 1, 1, 1, 4, 4, 4, 4, 7
\}~~\text{and}~~\operatorname{DS}(\Gamma(G))=\{ 1, 2, 2, 2, 4, 4, 4, 7\}.$$

In point of fact, $D_8$ and $G=\mathcal V_4\rtimes C_2$ are isomorphic.  To see this, we invoke~(23.4) on page~107 of~\cite{asch}, which states:

\begin{quote}
Let $G$ be a nonabelian group of order $p^n$ with a cyclic subgroup of
index~$p$.  Then $G\cong\text{Mod}_{p^n}$, $D_{2^n}$, $SD_{2^n}$, or
$Q_{2^n}$.
\end{quote}
The paragraph immediately preceding this result states the modular group
$\text{Mod}_8=D_8$, while two paragraphs further back states the
semidihedral group $SD_{2^n}$ is only defined for $n\geq 4$.  Now,
consider an arbitrary nonabelian group $H$ of order~$8$.  Since every
group of exponent~$2$ is necessarily abelian\footnote{Proof:  $xxyy =
x^2y^2 = 1 = (xy)^2 = xyxy$ implies $xy=yx$ over all $x,y$ in the group.}, $H$ must
have an element of order~$4$.  It therefore has a cyclic subgroup of
index~$2$ and thus satisfies the theorem.  Consequently, there are
exactly~$2$ nonabelian groups of order~$8$: the dihedral group $D_8$ and
the quaternion group $Q_8 = \{\pm 1,\pm i,\pm j,\pm k\}$.  The only
involution in $Q_8$ is $-1$.  Therefore, whenever a nonabelian group $H$ 
of order~$8$ has more than one element of order~$2$, then necessarily $H$ 
is dihedral.  The subgroup $\mathcal V_4$ of $G=\mathcal V_4\rtimes C_2$ 
already contains 3 involutions; $G$ is therefore not quaternion, 
and is then necessarily dihedral.

\section{The Errors}

Given that the same group yielded distinct graphs, where are the logical
flaws in~\cite{wrong}?  One place where fault can be found is in the
quoting of the Normal Form Theorem (NFT).  Although not mentioned in the
article~\cite{wrong}, the open letter makes explicit
that the NFT being referenced is found on page~31 of~\cite{cohen}\footnote
{A look at \cite{cohen}'s index reveals four NFTs; there was thus need for 
external verification as to which was being used.}. This NFT is in regards 
to an amalgamated free product, not to a semidirect product.  Another issue 
at point, which comes up in the given examples above, is that a semidirect 
factorization of groups is not unique; a group $G$ can split over nonisomorphic, 
normal subgroups~$N_1,N_2$ even with $|N_1|=|N_2|$.  In short, there is no clear way to adjust 
this graph's definition to allow the resulting object to reflect the 
underlying group's structure.

\end{document}